\documentclass{svmult-ddm}

\usepackage{makeidx}         
\usepackage{graphicx}        
\usepackage{multicol}        
\usepackage[bottom]{footmisc}
\usepackage{url}

\renewcommand{\PackageWarningNoLine}[2]{}
\usepackage{amsmath}
\usepackage{amsfonts}
\usepackage{float}

\newcommand{\pdl}[2]{\frac{\partial #1}{\partial #2}}

\newcommand{\checked}[1][\null]{\ensuremath{\boldsymbol{\surd}}}

\begin{document}

\title*{A stochastic domain decomposition method for time dependent mesh generation}
\titlerunning{Stochastic mesh generation}
\author{Alexander Bihlo\inst{1}\and
Ronald D. Haynes\inst{2}}
\institute{Memorial University of Newfoundland, St. John's, NL, Canada
\texttt{abihlo@mun.ca}
\and Memorial University of Newfoundland, St. John's, NL, Canada \texttt{rhaynes@mun.ca}}
%
%
\maketitle

\vskip 0.2in
\section{Introduction}

We are interested in PDE mesh generation where the mesh is computed by solving a mesh PDE which
is coupled to the physical PDE of interest.
In~\cite{bihl13a} we proposed a stochastic domain decomposition (DD) method to find adaptive meshes by solving a linear elliptic mesh
generator.  The stochastic DD approach, as originally formulated in~\cite{aceb05a}, relies on an accurate numerical solution using the probabilistic form of the exact solution of the linear elliptic boundary value problem.
Monte--Carlo simulations are used to evaluate this probabilistic form of the solution only at the sub-domain interfaces.  These interface approximations can be computed independently and  are then used as Dirichlet boundary conditions for the deterministic sub-domain solves.
Within the framework of grid adaptation it is generally not necessary to solve the mesh PDEs with
high accuracy. The reason is that the mesh equations are only a means to an end.  Only a good quality mesh, one that allows an accurate representation of the physical PDE, is required.
This relatively low accuracy requirement makes the proposed stochastic DD method computationally more attractive, reducing the number of Monte--Carlo simulations required.

 Grid adaptation by a stochastic DD approach does generate interesting issues in its own right.   Grid quality
should be monitored during the interface solves to give a suitable stopping criteria for the stochastic portion of the algorithm.   In \cite{bihl13a} only the steady grid generation problem was considered.  Of course, in practice, the
problem of grid generation is coupled with the process of solving the system of physical, usually time
dependent, PDEs.  It is this latter issue that we begin to explore in this paper.

We are interested in time dependent PDEs whose solutions evolve on disparate space and time scales.  The solution
behaviour lends itself to the use of time dependent meshes which automatically adapt and evolve to efficiently
resolve the solution features.   The generation of these time dependent grids can be done either by statically
applying an elliptic mesh generator
using the physical solution obtained at the previous time step or by employing
 a time relaxation of the static mesh PDE resulting in a parabolic mesh equation, as in \cite{huan10a}.

In \cite{bihl13a} we applied the stochastic DD algorithm to a linear elliptic mesh generator.  Here we consider the
extension to time dependent mesh PDEs.
This extension to (linear) parabolic mesh generators is possible due to the existence of a stochastic representation of the exact solution of such linear parabolic problems.
For the sake of illustration, we will work with the time-relaxed form of the Winslow--Crowley variable
diffusion mesh generation method, first described  in~\cite{wins66a}.

\section{Winslow's method}\label{sec:winslow}

The {\em equipotential} method of
mesh generation in 2D, \cite{crowley62}, found the mesh
lines in the physical co-ordinates $x$ and $y$
as the level curves of the potentials $\xi$ and
$\eta$ satisfying  Laplace's equations
\begin{equation}\label{eq:crowley}\nabla^2 \xi  =0, \qquad \nabla^2 \eta = 0,\end{equation}
and  appropriate boundary conditions which ensure grid lines lie along the boundary of
the domain.  Here derivatives are with respect to the physical co--ordinates.
The physical mesh transformation,
$x(\xi,\eta)$ and $y(\xi,\eta)$, in the physical domain
$\Omega_p$,
can be found by (inverse) interpolation of the solution of (\ref{eq:crowley}) onto
 a (say) uniform $(\xi,\eta)$ grid.
In practice, the inversion to the physical co--ordinates is not necessary. Instead one could transform the physical PDE of interest into the computational co--ordinate system.
Winslow \cite{wins81a} generalized (\ref{eq:crowley})
by adding a diffusion coefficient $1/w(x,y) > 0 $ depending on the
gradient or other aspects of the solution.  This gives
the linear elliptic mesh generator
\begin{equation}\label{eq:linearmesh}
-\nabla\cdot\left(\frac{1}{w}\nabla \xi\right) = 0\quad\text{and}\quad -\nabla\cdot\left(\frac{1}{w}\nabla \eta\right) = 0.
\end{equation}
The
function $w$, known as a mesh density function,  characterizes regions where additional
mesh resolution is needed and in general depends on the solution of the physical PDE.

Here we assume the solution of the physical PDE is time dependent and hence the mesh density function is
changing with time, $w=w(t,x,y)$.   One could still use (\ref{eq:linearmesh}) to solve the mesh
transformation at each time $t$.    For time dependent PDEs this would result
a system of differential--algebraic equations for the physical solution and the mesh.  Instead, we choose to
relax (\ref{eq:linearmesh}) to obtain a parabolic linear mesh
generator of the form
\begin{align}\label{eq:RelaxedWinslowSystem}
 \xi_t=-\frac{1}{w} \nabla w\cdot\nabla\xi+\nabla^2\xi\quad\text{and}\quad \eta_t = -\frac{1}{w}\nabla w\cdot\nabla \eta+\nabla^2\eta.
\end{align}
This gives a mesh PDE which depends explicitly on the mesh speed and provides a degree of temporal
smoothing for the mesh, cf.\ \cite{huan94a}.

\section{Linear parabolic differential equations and stochastic domain decomposition}

The system of mesh PDEs (\ref{eq:RelaxedWinslowSystem}) is of the form
\begin{equation}\label{eq:ParabolicMeshSystemGeneral}
 \xi_t=\mathrm{L}\xi,\quad \eta_t=\mathrm{L}\eta,
\end{equation}
where $\xi(t,x,y)$ and $\eta(t,x,y)$ are the computational coordinates defined over $[0,T]\times\Omega_{\rm p}$, where $\Omega_{\rm p}$ is the spatial domain in physical coordinates. In system~\eqref{eq:ParabolicMeshSystemGeneral}, $\mathrm{L}$ is a linear elliptic operator of the form
\[
 \mathrm{L}=a_{ij}\frac{\partial^2}{\partial x_i\partial x_j}+b_i\pdl{}{x_i},
\]
with continuous coefficient matrix $a(t,x,y)=(a_{ij})(t,x,y)$, $i,j\in\{1,2\}$, and drift vector $\mathbf{b}=(b_1,b_2)^{\rm T}(t,x,y)$. Here we employ the summation convention over repeated indices.

System~\eqref{eq:ParabolicMeshSystemGeneral} is accompanied by boundary and initial conditions
$
 \xi|_{\partial\Omega_{\rm p}}=f(t,x,y),\,\, \eta|_{\partial\Omega_{\rm p}}=g(t,x,y),\,\,
\xi(0,x,y)=\xi_0(x,y),\,\,\text{and}\,\, \eta(0,x,y)=\eta_0(x,y).
$

The solution of such linear parabolic problems can be described using the tools of stochastic calculus \cite{aceb10a,kara91a}.
The point--wise
solution of system~\eqref{eq:ParabolicMeshSystemGeneral} at $(t,x,y)\in[0,T]\times\Omega_p$ is given probabilistically as
\begin{align}\label{eq:ProbabilisticSolutionGeneral}
\begin{split}
 \xi(t,x,y)=\mathrm{E}\left[\xi_0(\mathbf{X}(t))\mathbf{1}_{[\tau_{\partial\Omega_{\rm p}}>t]} \right] + \mathrm{E}\left[f(t-\tau_{\partial\Omega_{\rm p}}, \mathbf{X}(\tau_{\partial\Omega_{\rm p}}))\mathbf{1}_{[\tau_{\partial\Omega_{\rm p}}<t]} \right],
\end{split}
\end{align}
where the process $\mathbf{X}(t)=(x(t),y(t))^{\rm T}$ satisfies, in the \^{I}to sense, the stochastic differential equation (SDE)
\[
 \mathrm{d}\mathbf{X}(t)=\mathbf{b}(t,\mathbf{X}(t))\mathrm{d}t+\sigma(t,\mathbf{X}(t))\mathrm{d}\mathbf{W}(t).
\]
The relation between $\sigma$ and $(a_{ij})$ is given through
\[
 \frac12\sigma(t,x,y)\sigma(t,x,y)^{\rm T}=a(t,x,y)
\]
for all $(t,x,y)\in[0,T]\times\mathbf{R}^{2}$. The solution for $\eta(t,x,y)$ is completely analogous.

In (\ref{eq:ParabolicMeshSystemGeneral}), the $\mathrm{E}[\cdot]$ denotes the expected value, $\tau_{\partial\Omega_{\rm p}}$ is the time when the stochastic path starting at $(x,y)$ first hits the boundary of the physical domain $\Omega_{\rm p}$, $\mathbf{W}$ is two-dimensional Brownian motion and $\mathbf{1}$ is the indicator function.

The time dependent mesh generator~\eqref{eq:RelaxedWinslowSystem} is a special case of the general form~\eqref{eq:ParabolicMeshSystemGeneral} with
\begin{equation}\label{eq:ParametersWinslowSystem}
 a(t,x,y)=I_2,\quad b_1(t,x,y)=-\frac1w w_x,\quad b_2(t,x,y)=-\frac1w w_y,
\end{equation}
where $I_2$ is the $2\times2$ identity matrix.

For our two dimensional mesh generator we choose the initial conditions $\xi(t=0,x,y)=\xi_0(x,y)=x$ and $\eta(t=0,x,y)=\eta_0(x,y)=y$, corresponding to an initial uniform mesh,  and the static boundary conditions $\xi(t,x_{\rm l},y)=0$, $\xi(t,x_{\rm r},y)=1$, $\eta(t,x,y_{\rm l})=0$ and $\eta(t,x,y_{\rm u})=1$.  This ensures we use the standard computational domain $\Omega_{\rm c}=[0,1]\times[0,1]$ and the rectangular physical domain $\Omega_{\rm p}=[x_{\rm l},x_{\rm r}]\times[y_{\rm l},y_{\rm u}]$.
The remaining boundary conditions for $\xi(t,x,y_{\rm l}),\xi(t,x,y_{\rm u}),\eta(t,x_{\rm l},y)$ and $\eta(t,x_{\rm r},y)$ are determined by solving the 1D version of (\ref{eq:linearmesh}) along the boundaries. Collectively, we use $f$ and $g$ to denote
these boundary conditions for $\xi$ and $\eta$.

Hence we have to solve the SDE
\begin{subequations}\label{eq:StochasticSolutionOneDimensional}
\begin{equation}\label{eq:StochasticSolutionOneDimensionalA}
 \mathrm{d}\mathbf{X}(t)=-\frac1w \nabla w\,\mathrm{d}t+\sqrt{2}\,\mathrm{d}\mathbf{W}(t),\quad
\end{equation}
for the single path $X(t)$.
The stochastic form of the exact solution of Eq.~\eqref{eq:RelaxedWinslowSystem} for $\xi $ is then obtained
by evaluating
\begin{equation}\label{eq:StochasticSolutionOneDimensionalB}
 \xi(t,x,y)=\mathrm{E}\left[\xi_0(\mathbf{X}(t))\mathbf{1}_{[\tau_{\partial\Omega_{\rm p}}>t]} \right] + \mathrm{E}\left[f(\mathbf{X}(\tau_{\partial\Omega_{\rm p}}))\mathbf{1}_{[\tau_{\partial\Omega_{\rm p}}<t]} \right].
\end{equation}
\end{subequations}
Thepoint--wise solution for $\eta(t,x,y)$ is obtained in an analogous fashion.

In principle, the probabilistic solution~\eqref{eq:StochasticSolutionOneDimensional} allows one to determine the computational coordinates $\xi$ and $\eta$ at each point in the space--time domain $[0,T]\times\Omega_{\rm p}$. However, this is prohibitively expensive. A more efficient approach is to evaluate the  solution~\eqref{eq:StochasticSolutionOneDimensional} only at certain points in space and time which then serve as boundary points for a
DD implementation. This stochastic DD approach for parabolic problems has been studied by Acebr{\'o}n et. al.~\cite{aceb10a}.




In the mesh generation context it is not possible to obtain the solution of~\eqref{eq:ProbabilisticSolutionGeneral} at all times, as the solution of the mesh PDE is coupled to the physical solution. That is, rather than solving~\eqref{eq:ProbabilisticSolutionGeneral} for a time $t\in[0,T]$, it is generally only be possible to use this stochastic solution to advance the solution of~\eqref{eq:ParabolicMeshSystemGeneral} over one single time step from $t^n$ to
$t^{n+1}$.
In this case, $\xi_0$ and $\eta_0$ should be interpreted as the values of $\xi$ and $\eta$ at time $t^n$
and the monitor function, $w$, is given at either $t^n$ or $t^{n+1}$ and remains constant over the time step.

\section{The numerical method}


\textbf{Stochastic solver and domain decomposition.} The use of the stochastic solution~\eqref{eq:ProbabilisticSolutionGeneral} for the time-relaxed Winslow mesh generator with parameters~\eqref{eq:ParametersWinslowSystem} is straightforward. We solve~\eqref{eq:StochasticSolutionOneDimensionalA} using the classical Euler--Maruyama scheme, i.e.\ we employ linear time-stepping. An alternative would be to use exponential time-stepping as advocated e.g.\ in~\cite{aceb05a,bihl13a,jans03a}.  In our results linear time-stepping gives sufficient accuracy.  The components of the Brownian motion $\mathrm{d}\mathbf{W}(t)$ are computed as $\sqrt{\Delta t}\,\mathcal{N}(0,1)$, where $\mathcal{N}(0,1)$ is a normally distributed random number with mean zero and variance one~\cite{kara91a}.

The time dependent weight only becomes available only at each time step (due to a possible coupling with a physical PDE).  Hence we are only able to employ formula~\eqref{eq:StochasticSolutionOneDimensionalB} to integrate over a single time step, i.e.\ from $t^n$ to $t^{n+1}$. Over this time step, the weight function is evaluated at $t^n$ and held constant, i.e.\ we have $w^n(x,y)=w(t^n,x,y)$ in~\eqref{eq:StochasticSolutionOneDimensionalA}. Accordingly, $\xi_0$ in Eq.~\eqref{eq:StochasticSolutionOneDimensionalB} is to be interpreted as $\xi_0^n=\xi(t^n,x,y)$, i.e.\ the values of the computational coordinates at the current time $t^n$.
Moreoever, the boundary functions $f$ and $g$ are updated at each time to reflect changes in the physical solution.

We then numerically integrate the SDE~\eqref{eq:StochasticSolutionOneDimensionalA} from $t^n$ to $t^{n+1}$. The drift vector $\mathbf{b}=-\nabla w/w$ is required everywhere along the path of the stochastic process $\mathbf{X}(t)$ but is only available at the grid points of the domain. Bilinear interpolation is used to obtain the values of $\mathbf{b}$ in between these grid points. The derivatives in $\nabla w$ are approximated using finite differences.

In the DD context, the stochastic solution is only required at a selection of points, $(x_k^i,y_k^i)$, which live on the
interfaces between sub-domains.
One time step $\Delta t$ is split into several smaller sub-time steps in order to numerically integrate the SDE~\eqref{eq:StochasticSolutionOneDimensionalA} from $t^n$ to $t^{n+1}$. We found this splitting of $\Delta t$ into sub-time steps necessary to determine sufficiently accurate whether the stochastic processes started at an
interface point
has left the domain $\Omega_{\rm p}$ during $\Delta t$.   This is not unlike the $M^k$ approach for mesh generation discussed in \cite{huan10a}.
At each sub-time step, a boundary test is performed to determine whether the stochastic process has left the domain $\Omega_{\rm p}$. If this is the case, the process contributes via the second term in Eq.~\eqref{eq:StochasticSolutionOneDimensionalB} to the approximation of $\xi(t^{n+1},x_k^i,y_k^i)$.
If the stochastic process did not leave the domain until $t^{n+1}$ is reached, it contributes to the first term in the approximation of $\xi(t^{n+1},x_k^i,y_k^i)$  in Eq.~\eqref{eq:StochasticSolutionOneDimensionalB}.
to the approximation of $\xi(t^{n+1},x_k^i,y_k^i)$.
The computation of $\eta(t^{n+1},x_k^i,y_k^i)$ is handled in the analogous way.
The expected values are then replaced by arithmetic means and approximated using the Monte-Carlo method.
Note, it is not desirable to make $\Delta t$ itself smaller, as this would degrade the efficiency of the (deterministic) implicit single domain solver, which is described below.

\medskip

\textbf{Deterministic single-domain solver.} The values of $\xi$ and $\eta$ along the subdomain interfaces
serve as boundary conditions for the single-domain solver. The single-domain solver we employ is an implicit finite-difference discretization of Eq.~\eqref{eq:RelaxedWinslowSystem}. The matrix system is solved using LU-factorization.

\medskip

\textbf{Parallelization and further speed-up.} It is well-known that Monte-Carlo techniques converge rather slowly~\cite{pres07a} and are usually most competitive for problems in high dimensions. The use of the stochastic solution to obtain the interface values of a DD problem only,  however, is considerably more efficient and provides a fully parallel grid generation algorithm. In particular, it is not required to pass information from one sub-domain to another. Moreover, the stochastic solutions on the interfaces can be determined at each point separately and each Monte-Carlo simulation is independent. Additionally, each sub-domain solution could potentially be assigned to a single processor once the interface solutions are obtained, yielding excellent scalability. Due to the fully parallel nature of the algorithm, the method is also fault tolerant. This renders the method suitable for an implementation on massively parallel computing architectures, cf.~\cite{aceb05a,aceb10a,bihl13a}.

A further source of improvement stems from the fact that not all values of $\xi$ and $\eta$ on the interfaces have to be computed using the stochastic solution~\eqref{eq:StochasticSolutionOneDimensional}. As proposed in~\cite{aceb05a} it may be sufficient to use the stochastic solution only at few points on the interface and recover the solution at the remaining interface points using interpolation. In~\cite{bihl13a} we have used a relatively simple optimal placement strategy to determine the most important locations on the interface where the stochastic solution should be computed. We use the same strategy in the present algorithm, i.e.\ the stochastic solution is computed near the maxima and minima of $\rho_x$ and $\rho_{xx}$ along the horizontal interfaces and $\rho_y$ and $\rho_{yy}$ along the vertical interfaces.

\section{Numerical Results}

\begin{figure}[H]\label{fig:somegrids}
\includegraphics[width=\textwidth,height=5in]{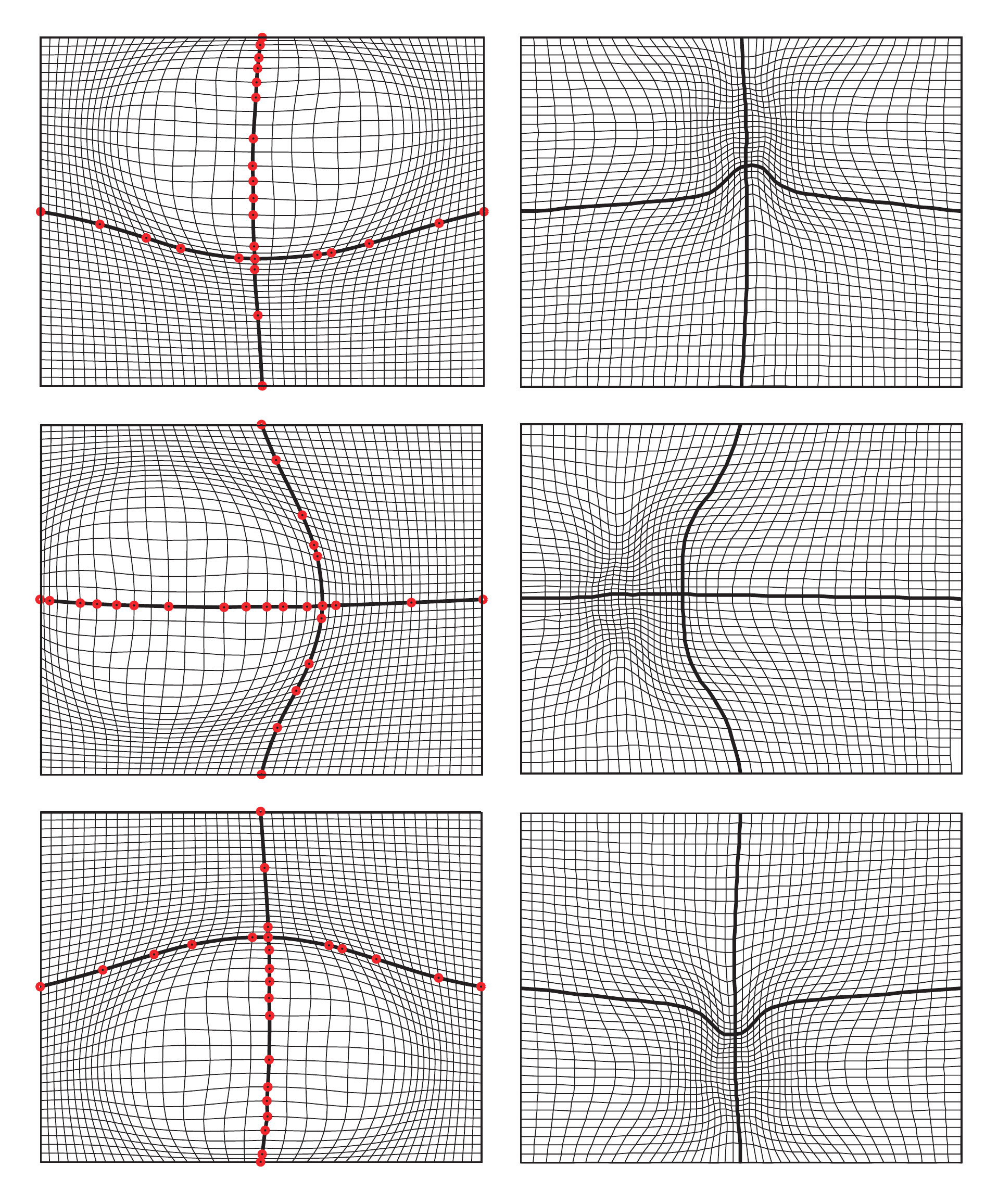}
\caption{Top to bottom: Meshes obtained from the parabolic mesh generator~\eqref{eq:RelaxedWinslowSystem} using the deterministic--stochastic method at $t=0.25$, $t=0.5$, and  $t=0.75$. Left: Meshes over the physical domain. Right: Meshes over the computational domain obtained from the former using natural neighbor interpolation. Thick line: Sub-domain interfaces. Circles: Points where the mesh is obtained using the stochastic solution~\eqref{eq:StochasticSolutionOneDimensional}.}
\end{figure}

We present an example our combined deterministic-stochastic DD method to generate an adaptive (moving) mesh
for the weight function $w=1/\rho$, where
\[
 \rho=1+\alpha\exp\left(\beta\left|\left(x-\frac12-\frac14\cos(2\pi t)\right)^2+\left(y-\frac12-\frac14\sin(2\pi t)\right)^2-\frac1{100}\right|\right).
\]
We choose the parameters $\alpha=10$ and $\beta=-50$ used in~\cite{huan10a}. Both the physical and computational domain are the unit square. The grid we generate has $41\times 41$ nodes and is divided into four sub-domains. On the interfaces we determine the stochastic solution at the key points using the optimal placement strategy mentioned in the previous section. Piecewise cubic Hermite interpolation is used to determine the remaining interface points. We integrate~\eqref{eq:RelaxedWinslowSystem} up to $t=0.75$ using $\Delta t=0.001$. Each time step is split into $20$ sub-time steps while solving the SDE~\eqref{eq:StochasticSolutionOneDimensionalA} and $N=10000$ Monte-Carlo simulations are used to estimate the expected value in~\eqref{eq:StochasticSolutionOneDimensionalB}. The resulting meshes at $t=0.25$, $t=0.5$, and $t=0.75$ are depicted in Fig.~\ref{fig:somegrids}.

The method is able to produce smooth meshes over the physical domain that adapt well to the time-dependent monitor function. No explicit smoothing was applied to the final meshes in this example.  In general we have found
sub-domain smoothing to be a way to further reduce the number of Monte-Carlo simulations needed in the
probabilistic expression~\eqref{eq:StochasticSolutionOneDimensionalB}, see~\cite{bihl13a}.

\section{Conclusion}

In this paper we have proposed a new deterministic--stochastic domain decomposition method for the construction of adaptive moving meshes suitable for time-dependent problems. The method is fully parallelizable as the values of the computational coordinates $\xi$ and $\eta$ on the single sub-domains can be determined without information exchange from neighboring sub-domains.

\looseness=-1
Future refinements include the use of exponential time-stepping 
to solve the SDE~\eqref{eq:StochasticSolutionOneDimensionalA}.  More generally,  more sophisticated boundary tests could better determine the first exit time of a stochastic process. This will allow using larger time steps in the solution of~\eqref{eq:StochasticSolutionOneDimensionalA} thus making the method more efficient.
An alternate approach to generate time dependent meshes is to apply the stochastic--DD method from \cite{bihl13a} to the sequence of elliptic problems which result from discretizing (\ref{eq:linearmesh}) in time.

{\bf Acknowledgements.} This research was supported by the Natural Sciences and Engineering Research Council of Canada (NSERC). The authors thank Professor Weizhang Huang (Kansas) for helpful remarks.

%
\bibliographystyle{spmpsci}
\bibliography{haynes_contrib_4} 



\end{document}